\documentclass[12pt]{amsart}

\usepackage{amsmath}
\usepackage{latexsym}
\usepackage{amssymb}
\usepackage{amsfonts}

\usepackage{graphics}

\renewcommand{\proof}{\par\noindent{\it Proof.\ \ }}
\def\qed{\ifmmode\square\else\nolinebreak\hfill
$\square$\fi\par\vskip5pt}

\def\l{\langle} \def\r{\rangle} 
\def\div{\,\big|\,}

 \def\ZZ{\mathbb Z}

\def\Aut{{\sf Aut}}  \def\Out{{\sf Out}}

\def\K{{\sf K}}
\bf

\def\D{{\rm D}}  
 \def\S{{\rm S}} 
\def\J{{\rm J}} \def\M{{\rm M}}
\def\soc{{\rm soc}} 
\def\C{{\bf C}}\def\N{{\bf N}}\def\Z{{\bf Z}} \def\O{{\bf O}}

\def\Ga{{\it\Gamma}}  

\def\a{\alpha}

\def\GammaL{{\rm \Gamma L}}

\def\PGammaL{{\rm P\Gamma L}} 
\def\A{{\rm A}}\def\Sym{{\rm Sym}}

\def\PSL{{\rm PSL}}\def\PGL{{\rm PGL}}
\def\GL{{\rm GL}}

 \def\PSU{{\rm PSU}} 

\def\Sz{{\rm Sz}} \def\McL{{\rm McL}} 
\def\Ru{{\rm Ru}} 
  \def\D{{\rm D}}

\def\HS{{\rm HS}}

\def\ON{{\rm O'N}}

\newtheorem{theorem}{Theorem}[section]%
\newtheorem{lemma}[theorem]{Lemma}%
\newtheorem{corollary}[theorem]{Corollary}%
\newtheorem{example}[theorem]{Example}%

\begin{document}

\title[edge-primitivity]{On edge-primitive and $2$-arc-transitive graphs}
\thanks{Supported by the National Natural Science Foundation of China (11731002) and the Fundamental Research Funds for the Central Universities.}
\author[Lu]{Zai Ping Lu}

\address{
Zaiping Lu\\
Center for Combinatorics, LPMC,
Nankai University,
Tianjin 300071, China
}
\email{lu@nankai.edu.cn}

\begin{abstract}

A graph is edge-primitive if its automorphism group acts primitively on the edge set. In this short paper, we prove that a finite $2$-arc-transitive edge-primitive graph has almost simple automorphism group if it is neither a cycle nor a complete bipartite graph. We also present two examples of such graphs, which are $3$-arc-transitive and have faithful vertex-stabilizers.

\vskip 10pt

\end{abstract}
\maketitle
\date\today

\vskip 30pt

\section{Introduction}

All graphs and groups considered in this paper are assumed to be finite.

\vskip 5pt

A graph in this paper is a pair $\Ga=(V,E)$
of a nonempty set $V$ and a set $E$ of $2$-subsets of $V$. The elements in $V$ and $E$ are called the vertices and edges of $\Ga$, respectively. The number $|V|$  of vertices is called the order of $\Ga$. For $v\in V$, the set $\Ga(v)=\{u\in V\mid \{u,v\}\in E\}$ is  called the neighborhood of $v$ in $\Ga$, while $|\Ga(v)|$ is the valency of $v$. We say that $\Ga$ has valency $d$ or $\Ga$ is $d$-regular if its vertices all have equal valency $d$. For an integer $s\ge 1$, an $s$-arc in $\Ga$ is
an $(s+1)$-tuple $(v_0,v_1,\ldots, v_s)$ of vertices such that $\{v_i,v_{i+1}\}\in E$ and $v_i\ne v_{i+2}$ for all possible $i$. A $1$-arc is also called an arc.

\vskip 5pt

Let $\Ga=(V,E)$ be a graph. A permutation $g$ on $V$ is called  an automorphism of $\Ga$ if $\{u^g,v^g\}\in E$ for all $\{u,v\}\in E$. Let $\Aut\Ga$ denote the set of all automorphisms of $\Ga$.
Then $\Aut\Ga$ is a subgroup of the symmetric group $\Sym(V)$, and called the automorphism group of $\Ga$. Note that the group $\Aut\Ga$ has a natural action on the edge set $E$ (and also on the set of $s$-arcs).
The graph $\Ga$ is called {\em edge-transitive} if $E\ne\emptyset$ and for each pair of  edges there
exists some $g\in \Aut\Ga$ mapping one of these two edges to the other one.
(Similarly, we may define {\em vertex-transitive}, {\em arc-transitive} or {\em $s$-arc-transitive} graphs.)
An edge-transitive graph is called {\em edge-primitive} if  some (and hence every) {\em edge-stabilizer}, the subgroup of its automorphism group which fixes a given edge, is a maximal subgroup of the  automorphism group.

\vskip 5pt

It is well-known that
edge-transitive graphs and hence edge-primitive graphs are either bipartite or vertex-transitive. As a subclass of the edge-transitive graphs, edge-primitive graphs posses
more restrictions on their symmetries and  automorphism groups.
For example, a connected edge-primitive graph is necessarily   arc-transitive
provided that it is not a star graph. In \cite{Giu-Li-edge-prim}, appealing
to the O'Nan-Scott Theorem for (quasi)primitive groups \cite{Prag-o'Nan}, Giudici and Li investigated the structural properties of  edge-primitive graphs, particularly, on their automorphism groups. Let $\Ga=(V,E)$ be an arc-transitive and edge-primitive graph which is neither a cycle nor a complete bipartite graph.
If $\Ga$ is bipartite then let $\Aut^+\Ga$ be the subgroup of $\Aut\Ga$ preserving the bipartition.
By \cite{Giu-Li-edge-prim}, as a primitive group on $E$, only $4$ of the eight O'Nan-Scott types for (quasi)primitive groups may occur  for $\Aut\Ga$, say {\rm SD}, {\rm CD}, {\rm PA} and {\rm AS}. For the first two types, $\Ga$ is bipartite and $\Aut^+\Ga$  is quasiprimitive of type {\rm CD} on each bipartite half. For the last two types, with one exception case, $\Aut\Ga$ or $\Aut^+\Ga$
is quasiprimitive on $V$ or on  each bipartite half respectively of the same type for $\Aut\Ga$ on $E$. In this paper, we will work on the types of  $\Aut\Ga$ on $E$ and on $V$ under the further assumption that $\Ga$ is $2$-arc-transitive.

\vskip 5pt

The interests for edge-primitive graphs arises partially from the fact that
many (almost) simple groups may be represented as  the automorphism groups of
edge-primitive graphs. Consulting the Atlas \cite{Atlas}, one may get first-hand such examples.
For example, the sporadic Higman-Sims group $\HS$ is the automorphism group of
a rank $3$ graph with order $100$ and valency $22$, which is in fact a $2$-arc-transitive and edge-primitive graph; the sporadic Rudvalis group $\Ru$
is the automorphism group of a rank $3$ graph with order $4060$ and valency $2304$, which is  edge-primitive  but not $2$-arc-transitive. Besides, the almost groups $\PSU(3,5).2$, $\M_{22}.2$, $\J_2.2$ and $\McL.2$ all have representations on edge-primitive graphs. The reader may refer to \cite{5val, 4val, s-arc, pval, 3val} for more examples of edge-primitive graphs which have almost simple automorphism groups.
Of course, using the constructions given in \cite{Giu-Li-edge-prim}, one can easily  construct examples of  edge-primitive graphs with automorphism groups not almost simple.

\vskip 5pt

We have a strong impression from the known examples for edge-primitive graphs in the literature that a $2$-arc-transitive and edge-primitive graph has almost simple automorphism group unless it is a cycle or a complete bipartite graph. Yet could it be so? Yes, it is true! We shall prove the following result in Section \ref{proof-sect}.

\begin{theorem}\label{main-result}
Let $\Ga=(V,E)$ be an edge-primitive $d$-regular graph for some $d\ge 3$.
If   $\Ga$ is $2$-arc-transitive, then either $\Ga$ is a complete bipartite graph, or $\Ga$ has almost simple automorphism group.
\end{theorem}

\noindent{\bf Remarks on Theorem \ref{main-result}}.

(1) Li and Zhang \cite{s-arc} proved that $4$-arc-transitive and edge-primitive graphs have almost automorphism groups. Further, as a sequence of their classification on almost simple primitive groups with soluble point-stabilizers, they give a complete list for $4$-arc-transitive and edge-primitive graphs.

(2) By Theorem \ref{main-result}, appealing to  the classification of almost simple groups with soluble maximal subgroups, it might be feasible to classify  $2$-arc-transitive and edge-primitive graphs with soluble edge-stabilizers.

\vskip 30pt

\section{Preliminaries}\label{pre}

For the subgroups of (almost) simple groups, we sometimes follow the notation used in the Atlas~\cite{Atlas}, while we also use $\ZZ_l$ and $\ZZ_p^k$ to denote respectively
the cyclic group of order $l$ and the elementary abelian group of order $p^k$.

\subsection{Primitive groups}
In this subsection, $\Omega$ is nonempty finite  set, and $G$   is a transitive  subgroup of the symmetric group $\Sym(\Omega)$. Let $\soc(G)$ be the socle of $G$,
that is, $\soc(G)$ is generated by all minimal normal subgroups of $G$.

\vskip 5pt

Consider the point-stabilizer $G_\a:=\{g\in G\mid \a^g=\a\}$, where $\a\in \Omega$. Then
\begin{enumerate}
\item[(1)] $G$ is primitive if $G_\a$ is a maximal subgroup of $G$;
\item[(2)] $G$ is ${3\over 2}$-transitive if $G_\a$ is ${1\over 2}$-transitive on $\Omega\setminus\{\a\}$, that is, all  $G_\a$-orbits on $\Omega\setminus\{\a\}$  have equal length $>1$;
\item[(3)] $G$ is a Frobenius group if  $G_\a$ is semiregular on $\Omega\setminus\{\a\}$;
\item[(4)] $G$ is $2$-transitive if  $G_\a$ is transitive on $\Omega\setminus\{\a\}$.
\end{enumerate}
Note that (4) implies (1) and (2), and (2) implies (1) or (3) (refer to \cite[Theorem 10.4]{Wielandt}).

\vskip 5pt

Let $1\ne N\unlhd G$, a normal subgroup of $G$. Then $N$ is ${1\over 2}$-transitive, and $N_\a=N\cap G_\a\unlhd G_\a$, and so
 $G_\a$ is contained in the normalizer $\N_G(N_\a)$
  of $N_\a$ in $G$.
 Thus, if $G_\a$ is maximal then
 either $N_\a\unlhd G$  or $\N_G(N_\a)=G_\a$. The former case yields $N_\a=1$, while the latter case
 gives \[\N_N(N_\a)=N\cap \N_G(N_\a)=N\cap G_\a=N_\a.\]
  Then we have following simple fact for primitive groups.

\begin{lemma}\label{nor-prim}
If $G$ is primitive and $N\ne 1$ then  either $N$ is regular on $\Omega$ or $N_\a$ is self-normalized; if $G$ is $2$-transitive and $N\ne 1$
then  $N$ is either regular or  ${3\over 2}$-transitive on $\Omega$.
\end{lemma}

For an almost simple $2$-transitive group $G$, each non-trivial normal subgroup $N$ of $G$ is primitive, and in fact $2$-transitive except for the case where $N=\soc(G)=\PSL(2,8)$ acting on $28$ pionts, refer to
\cite[page 197, Table 7.4]{Cameron}. Next we consider the normal subgroups of affine $2$-transitive groups.
Refer to  \cite[page 195, Table 7.3]{Cameron} for a complete list of affine $2$-transitive groups.
We consider the affine $2$-transitive groups in their natural actions.

\begin{lemma}\label{nor-sub-affine}
Let $G$ be an affine  $2$-transitive group  and $1\ne N\unlhd G$.
 If $N$ is imprimitive on $\Omega$, then
$N$ is a soluble Frobenius group,
$N_0$ is cyclic, and either $G_0\le \GammaL(1,q)$  or  $N_0\le \Z(G_0)$, where $q$ is not a prime.
\end{lemma}
\proof
Assume that $N$ is imprimitive. Then $N\ne G$, and so $N_0\ne G_0$.
Further, by Lemma \ref{nor-prim} and \cite[Theorem 10.4]{Wielandt}, $N$ is a Frobenius group.
Let $|\Omega|=p^k$ for a prime $p$. We may write $G_0\le \GL(k,p)$, $G=\ZZ_p^k{:}G_0$ and $N=\ZZ_p^k{:}N_0$.
Since $N$ is imprimitive, $N_0$ is not maximal in $N$, and thus $N_0$ is a  normal reducible subgroup of $G_0$.
Then, by \cite[Lemma 5.1]{Hering}, $N_0$ is cyclic and $|N_0|$ is a divisor of $p^l-1$, where $l<k$ and $l\div k$.
Finally, the lemma follows from checking all affine $2$-transitive groups one by one.
\qed

If every minimal normal subgroup of  $G$ is transitive on $\Omega$, then $G$ called a quasiprimitive group.
Praeger \cite{Prag-o'Nan, Prag-o'Nan-survey} generalized the O'Nan-Scott Theorem for primitive groups to quasiprimitive
groups, which says that a  quasiprimitive
group has one of the following eight types: HA, HS, HC, TW, AS, SD, CD and PA. In particular,
if $G$ is quasiprimitive then $G$ has at most two minimal normal subgroups, and  if two (for HS and HC) then they are  isomorphic and regular.

\vskip 5pt

Suppose that $G$ has a transitive insoluble minimal normal subgroup $N$. Then $G=NG_\a$ for $\a\in \Omega$. Write $N=T_1\times \cdots \times T_k$ for isomorphic nonabelian simple groups $T_i$ and integer $k\ge 1$. Then $G_\a$ acts transitively on $\{T_i\mid 1\le i\le k\}$ by conjugation. Note that, for $g\in G_\a$ and $1\le i\le k$,
\[((T_i)_\a)^g=(T_i\cap G_\a)^g=T_i^g\cap G_\a^g=(T_i)^g_\a=(T_j)_\a \mbox{ for some } j.\]
Then $G_\a$ acts transitively on $\{(T_i)_\a\mid 1\le i\le k\}$ by conjugation.
Clearly, $(T_1)_\a\times \cdots \times (T_k)_\a\le N_\a$; however, the equality is not necessarily holds even if $G$ is quasiprimitive. A sufficient condition for this equality is that $G$ is primitive and of type AS or PA,
refer to \cite[Theorem 4.6]{Dixon} and its proof.
In survey, we have the simple fact as follows.

\begin{lemma}\label{PA-prim}
Assume that $G$ has a transitive   minimal normal subgroup $N=T_1\times \cdots\times T_k$, where $T_i$ are isomorphic nonabelian simple groups. Let $\a\in \Omega$.
Then $G_\a$ acts transitively on $\{(T_i)_\a\mid 1\le i\le k\}$ by conjugation.
If further $G$ is primitive and of type {\rm AS} or {\rm PA}, then $N_\a=(T_1)_\a\times \cdots\times (T_k)_\a$.
\end{lemma}

\vskip 5pt

\subsection{Locally-primitive graphs}
In this subsection, $\Ga$ is a connected $d$-regular graph for some $d\ge 3$, and $G\le \Aut\Ga$.
 Assume further that the graph $\Ga$ is $G$-locally primitive, that is,
$G_v$ acts primitively  on $\Ga(v)$ for all $v\in V$.

\vskip 5pt

Fix an edge $\{u,v\}\in E$.
Note that
$G_v$ induces a primitive permutation group $G_v^{\Ga(v)}$ (on $\Ga(v)$).
Let $G_v^{[1]}$ be the kernel of $G_v$ acting on $\Ga(v)$.
Then $G_v^{\Ga(v)}\cong G_v/G_v^{[1]}$.
Set $G_{uv}^{[1]}=G_u^{[1]}\cap G_v^{[1]}$. Then $G_v^{[1]}$ induces a normal subgroup of
$(G_{u}^{\Ga(u)})_v$ with the kernel $G_{uv}^{[1]}$.
Writing $G_v^{[1]}$ and $G_v$ in group extensions,
\[\tag{$\divideontimes$}  G_v^{[1]}=G_{uv}^{[1]}.(G_v^{[1]})^{\Ga(u)}, \,\,\, G_v=(G_{uv}^{[1]}.(G_v^{[1]})^{\Ga(u)}).G_v^{\Ga(v)},\,\, \, G_{uv}=G_v^{[1]}.(G_v^{\Ga(v)})_u.\]

\vskip 5pt

Assume that $G$ is transitive on $V$. Then $G_{uv}^{[1]}$ is a $p$-group for some prime $p$, refer to \cite{Gardiner-73}.
Note that $G$ is transitive on the arcs of $\Ga$. There is some element in $G$ interchanging $u$ and $v$.
This implies that $(G_v^{[1]})^{\Ga(u)}\unlhd (G_v^{\Ga(v)})_u\cong (G_{u}^{\Ga(u)})_v$. Thus  we have the following lemma.

\begin{lemma} Assume that $G$ is transitive on $V$, and $\{u,v\}\in E$. Then $G_{uv}^{[1]}$ is a $p$-group, and  $(G_v^{[1]})^{\Ga(u)}$
is isomorphic to a normal subgroup of a point-stabilizer in $G_v^{\Ga(v)}$. In particular, $G_v$ is soluble if and only if  $G_v^{\Ga(v)}$ is soluble.
\end{lemma}

The graph $\Ga=(V,E)$ is said to be $(G,s)$-arc-transitive if $\Ga$ has an $s$-arc and $G$ acts transitively on the set of $s$-arcs of $\Ga$, where $s\ge 1$.
Note that $\Ga$ is $(G,2)$-arc-transitive if and only if $G$ is transitive on $V$, and $G_v^{\Ga(v)}$ is a $2$-transitive group for some (and hence every) $v\in V$.
By \cite{Gardiner74,Weiss81a,Weiss81b}, we have the following result.

\begin{theorem}
\label{double-star}
Assume that $\Ga=(V,E)$ is $(G,2)$-arc-transitive. Then $\Ga$ is not $(G,8)$-arc-transitive. Further,
\begin{enumerate}
\item[(1)] if $G_{uv}^{[1]}=1$ then $\Ga$ is not $(G,4)$-arc-transitive.
\item[(2)] if $G_{uv}^{[1]}\ne 1$ then
 $G_{uv}^{[1]}$ is a nontrivial $p$-group,  $\O_p(G_{uv}^{\Ga(v)})\not=1$,
$\PSL(n,q)\unlhd G_v^{\Ga(v)}$, and $|\Ga(v)|={q^n-1\over q-1}$, where $n\ge 2$ and $q$ is a power of $p$; in this case,
$\Ga$ is $(G,4)$-arc-transitive if and only if $n=2$.
\end{enumerate}
\end{theorem}

\vskip 30pt

\section{The proof of Theorem \ref{main-result}}\label{proof-sect}
In this section, we let $\Ga=(V,E)$ be a connected graph of valency $d\ge 3$, and $G\le \Aut\Ga$. Assume that $\Ga$ is $G$-edge-primitive, that is, $G$ act primitively on $E$.
Then, by \cite[Lemma 3.4]{Giu-Li-edge-prim}, $G$ acts transitively on the arc set of $\Ga$.
Thus, for an edge $\{u,v\}\in E$, $d=|G_v:G_{uv}|$ and  $|G_{\{u,v\}}:G_{uv}|=2$.

\vskip 5pt

Let $1\ne N\unlhd G$. Then $N$ is transitive on $E$, and so either $N$ is transitive on $V$ or $N$ has  two orbits on $V$; for the latter case,
$N_v$ is transitive on $\Ga(v)$.
This implies that either $G=NG_v$, or   $|G:(NG_v)|=2$ and $N_{uv}=N_{\{u,v\}}$.
Note that $G=NG_{\{u,v\}}$ by the maximality of $G_{\{u,v\}}$ or the transitivity of $N$ on $E$. We have
\[
\begin{array}{lcl}
|G|&=&{|N||G_{\{u,v\}}|\over |N\cap G_{\{u,v\}}|}={|N||G_{\{u,v\}}|\over |N_{\{u,v\}}|}={2|N||G_{uv}|\over |N_{\{u,v\}}|}
={2|N||G_{v}|\over d|N_{\{u,v\}}|}\\
&=&{|N||G_{v}|\over |N_v|}\cdot {2|N_v|\over d|N_{\{u,v\}}|}=|NG_v|{2|N_v|\over d|N_{\{u,v\}}|}.
\end{array}
\]
Then the next lemma follows.
\begin{lemma}\label{nor-sub}
Let $1\ne N\unlhd G$. If $N$ is transitive on $V$ then $2|N_v|=d|N_{\{u,v\}}|$; if $N$ is intransitive on $V$ then $|N_v|=d|N_{\{u,v\}}|=d|N_{uv}|$. In particular, $N_v\ne 1$ and $N_v\ne N_{\{u,v\}}$.
\end{lemma}

Let $\K_{d,d}$ and $\K_{d+1}$ be the complete bipartite graph and complete   graph of valency $d$, respectively.

\begin{corollary}\label{M-v>d}
Let $1\ne N\unlhd G$.
Then either $\Ga\cong \K_{d,d}$,  or
 $N_{uv}\ne 1$
  and $N_{\{u,v\}}$ is self-normalized in $N$,  where $\{u,v\}\in E$.
\end{corollary}
\proof
Assume that $\Ga\not\cong \K_{d,d}$. Then, by the O'Nan-Scott Theorem and \cite[Lemmas 6.1, 6.2 and Propersition 6.13]{Giu-Li-edge-prim},
$G$ has no  normal subgroup acting regularly on $E$. Thus $N_{\{u,v\}}\ne 1$, and so $N_{\{u,v\}}$ is self-normalized in $N$
by Lemma \ref{nor-prim}.

Suppose that  $N_{uv}=1$. Then $N_{\{u,v\}}$ has order $2$, and so
$N_{\{u,v\}}\le \C_N(N_{\{u,v\}})\le \N_N(N_{\{u,v\}})=N_{\{u,v\}}$.
This implies that $\C_N(N_{\{u,v\}})=\N_N(N_{\{u,v\}})$, and then $N_{\{u,v\}}$ is a Sylow $2$-subgroup of $N$.
By   Burnside's transfer theorem (refer \cite[IV.2.6]{Huppert}), $N$ has normal $2'$-Hall subgroup, say $M$. Then this $M$ is normal in $G$ and regular on $E$, a contradiction.
\qed

\vskip 5pt

By \cite{Giu-Li-edge-prim}, if $\Ga\not\cong \K_{d,d}$ then
$G$ has type  {\rm SD}, {\rm CD}, {\rm AS} or {\rm PA} on $E$; in particular, $G$ has a unique (of course, insoluble) minimal normal subgroup. Thus,
  if $\Ga\not\cong \K_{d,d}$ then $G$ is insoluble, and so $G_{\{u,v\}}$ is not abelian by \cite[IV.7.4]{Huppert}. If  $G_{uv}$ is abelian the following result
  says that $\Ga\cong \K_{d,d}$ or $\K_{d+1}$.

\begin{theorem}\label{arc-stabilizer}
Assume that $\Ga\not\cong \K_{d,d}$.  Let $1\ne N\unlhd G$.
\begin{enumerate}
\item[(1)] If $N_{\{u,v\}}$
has a normal Sylow subgroup $P\ne 1$ then $P$ is also a Sylow subgroup of $N$; in particular,  $N_{\{u,v\}}$ is not abelain.
\item[(2)] If $N_{uv}$ is abelian then $N$ is transitive on the arc set of $\Ga$.
\item[(3)] If $N_{uv}$ is an abelian $2$-group then $\soc(G)=\PSL(2,q)$ and
$\Ga\cong \K_{q+1}$, where $q$ is a power of some prime with $q-1$  a power of $2$ greater than $8$.

\item[(4)] If $G_{uv}$ is an abelian group then $d=q$  and either $\soc(G)\cong\PSL(2,q)$  and
$\Ga\cong \K_{q+1}$, or $\soc(G)=\Sz(q)$, $\Aut\Ga=\Aut(\Sz(q))$ and $\Ga$ is $(\Sz(q),2)$-arc-transitive, where $q$ is a power of some prime.
\end{enumerate}
\end{theorem}
\proof
(1) Assume that $P\ne 1$ is a  normal Sylow $p$-subgroup of $N_{\{u,v\}}$.
Then $P$ is a characteristic subgroup of $N_{\{u,v\}}$, and so
$P\unlhd G_{\{u,v\}}$ as $N_{\{u,v\}}\unlhd G_{\{u,v\}}$. Thus $\N_G(P)\ge G_{\{u,v\}}$, and then $\N_G(P)=G_{\{u,v\}}$ by
the maximality of $G_{\{u,v\}}$. This gives $\N_N(P)=N\cap \N_G(P)=N\cap G_{\{u,v\}}=N_{\{u,v\}}$.
Choose a  Sylow $p$-subgroup $Q$ of $N$ with $P\le Q$. Then $\N_Q(P)\le Q\cap \N_G(P)=Q\cap N_{\{u,v\}}=P$. This yields $P=Q$, so   $P$ is a Sylow $p$-subgroup of $N$.

Suppose that  $N_{\{u,v\}}$ is abelian. Then  $N_{\{u,v\}}\le \C_N(P)\le \N_N(P)=N_{\{u,v\}}$, yielding $\C_N(P)=\N_N(P)$. By Burnside's transfer theorem,
$P$ has a normal complement $H$ in $N$, that is $N=PH$ with $P\cap H=1$ and $H\unlhd N$. Note that $H$ is a Hall subgroup of $N$. It follows that $H$ is characteristic in $N$, and hence $H\unlhd G$. Let $P$ runs over the Sylow subgroup of $N_{\{u,v\}}$. Then the resulting normal complements intersect at a normal complement of $N_{\{u,v\}}$ in $N$, which is normal in $G$ and regular on $E$.
This contradicts Corollary \ref{M-v>d}. Therefore, $N_{\{u,v\}}$ is nonabelian, and
(1) of this theorem follows.

\vskip 5pt

(2) Assume that $N_{uv}$ is abelian.
Then $N_{uv}\ne N_{\{u,v\}}$ by (1), and thus $(u,v)=(v,u)^x$ for some $X\in N_{\{u,v\}}$. Since $\Ga$ is $N$-edge-transitive, $\Ga$ is $N$-arc-transitive.

\vskip 5pt

(3) Assume that $N_{uv}$ is an abelian $2$-group.
Recall that $G$ has a unique minimal normal subgroup, say $M$. Then $M\le N$, and (1) and (2) hold for $M$. Then, since $M_{uv}$ is an abelian $2$-group,
$M_{\{u,v\}}$ is a Sylow $2$-subgroup of $M$, and
$M_{\{u,v\}}$ is not abelian.

Write $M=T_1\times \cdots\times T_k$, where $T_i$ are   isomorphic nonabelian simple groups. Recall that $M_{\{u,v\}}$ is a Sylow $2$-subgroup of $M$. For each $i$, choose a Sylow $2$-subgroup of $T_i$ with $Q_i\le M_{\{u,v\}}$.
Then $M_{\{u,v\}}=Q_1\times \cdots \times Q_k$. Noting that $Q_i$ are all isomorphic,
every $Q_i$ is nonabelian; otherwise, $M_{\{u,v\}}$ is abelian, a contradiction.
In particular, $Q_1\not\le M_{uv}$. Then
$M_{\{u,v\}}=M_{uv}Q_1$, and so \[Q_2\times \cdots \times Q_k\cong M_{\{u,v\}}/Q_1=M_{uv}Q_1/Q_1\cong M_{uv}/(M_{uv}\cap Q_1).\] Since $\M_{uv}$ is abelian, the only possibility is $k=1$. Thus $M=\soc(G)$ is simple.

By \cite[Corollary 5]{Goldschmidt},
$M_{\{u,v\}}$ has cyclic commutator subgroup.
Since $M_{\{u,v\}}$ is nonabelian, by \cite{Chabot}, $M$ is isomorphic to one of the Mathieu group $\M_{11}$,
$\PSL(2,q)$ (with $q^2-1$ divisible by $16$), $\PSL(3,q)$ (with $q$ odd) and $\PSU(3,q)$ (with $q$ odd).
If $M\cong \M_{11}$, then $G=M$, and so $M_{\{u,v\}}$ is maximal in $M$; however, by the Atlas \cite{Atlas}, a Sylow $2$-subgroup of $\M_{11}$ is not a maximal subgroup, a contradiction. Thus we next let $M\cong\PSL(2,q)$, $\PSL(3,q)$ or $\PSU(3,q)$.

Since $M$ is transitive on $E$, we know that $|E|=|M:M_{\{u,v\}}|$ is odd.
Thus $G$ is an almost simple primitive group (on $E$) of odd degree. Noting that $M_{\{u,v\}}=M\cap G_{\{u,v\}}$, by \cite{LS},
$M_{\{u,v\}}$ is known. Notice that the isomorphisms among simple groups (refer to \cite[Proposition 2.9.1 and Theorem 5.1.1]{KL}). Since $M_{\{u,v\}}$ is a Sylow $2$-subgroup of $M$,
the only possibility is that $M\cong \PSL(2,q)$, and
$M_{\{u,v\}}$ is
the stabilizer of some orthogonal  decomposition of a natural projective module
associated with $M$ into $1$-dimensional subspaces.
It follows that $M_{\{u,v\}}\cong \D_{q-1}$ or $\D_{q+1}$, and so $M_{uv}\cong \ZZ_{q-1\over 2}$ or $\ZZ_{q+1\over 2}$, respectively.
Since $M$ is transitive on the arcs of $\Ga$, we have $|M_v:M_{uv}|=d\ge 3$.
Checking the subgroups of $\PSL(2,q)$ (refer to \cite[II.8.27]{Huppert}),
we conclude that $M_{uv}\cong \ZZ_{q-1\over 2}$, $d=q$, $V=|M:M_v|=q+1$ and $M$ is $2$-transitive on $V$. Thus $\Ga\cong \K_{q+1}$.

\vskip 5pt

(4) Assume that $G_{uv}$ is abelian.
Let $M$ be the unique minimal normal subgroup of $G$. If  $M_{uv}$ is a $2$-group, then (4) of this theorem follows from (3).

We next assume that $|M_{uv}|$ has an odd prime divisor $p$.
By (1), the unique Sylow $p$-subgroup of $M_{uv}$ is also a Sylow $p$-subgroup of $M$.
Write $M=T_1\times \cdots\times T_k$ , where $T_i$ are  isomorphic nonabelian simple groups.
By (1) of this theorem, $\M_{\{u,v\}}$ is not abelian, so $M_{\{u,v\}}\not\le G_{uv}$, and then $G_{\{u,v\}}=M_{\{u,v\}}G_{uv}$. Thus $G=MG_{uv}$, and hence $G_{uv}$ acts transitively on $\{T_1,\ldots,T_k\}$ by conjugation. Choose, for each $i$,  a Sylow $p$-subgroup $P_i$ of $T_i$ such that $P_1\times\cdots\times P_k$ is the unique Sylow subgroup  $M_{uv}$. Since $G_{uv}$ is abelian,  we have $P_1=P_1^x\le T_1^x$ for $x\in G_{uv}$. It follows that $P_1\le T_i$ for all $i$.
The only possibility is that $k=1$, and so  $M$ is simple.

Note that $G$ is an almost simple  group with a soluble maximal subgroup $G_{\{u,v\}}$. Then, by \cite{s-arc}, both $M=\soc(G)$ and $M_{\{u,v\}}=M\cap G_{\{u,v\}}$  are  known.
Since $M_{\{u,v\}}$ has an abelian subgroup of index $2$, it follows that
either
$M\cong \PSL(2,q)$ and $M_{\{u,v\}}\cong \D_{2(q\pm 1)\over (2,q-1)}$,
 or $M=\Sz(q)$ and $M_{\{u,v\}}\cong \D_{2(q-1)}$. 
Check the subgroups of $M$, refer to \cite{Suzuki} for $\Sz(q)$.
The former case yields that $M_v\cong [q]{:}\ZZ_{q-1\over (2,q-1)}$ and $\Ga\cong \K_{q+1}$.
Assume that  $M=\Sz(q)$ and $M_{\{u,v\}}\cong \D_{2(q-1)}$. Then $M_v\cong [q]{:}\ZZ_{q-1}$ and $d=q$; in this case, $\Ga$ is $(M,2)$-arc-transitive.
By \cite{FP1}, we have that $\Aut\Ga=\Aut(\Sz(q))$ and $\Ga$ is unique up to isomorphism.
Thus 
(4) of this theorem follows.
\qed


\begin{lemma}\label{G-E-PA}
Assume  that $G$ has type {\rm PA} on $E$. Let $\soc(G)=T_1\times \cdots\times T_k$.
Then $(T_i)_{uv}\ne 1$ for each $i$ and $\{u,v\}\in E$; in particular, every $T_i$ is not semiregular.
\end{lemma}
\proof
Let $M=\soc(G)$.
By Lemma \ref{PA-prim}, $M_{\{u,v\}}=(T_1)_{\{u,v\}}\times \cdots\times (T_k)_{\{u,v\}}$, and $(T_i)_{\{u,v\}}$  all have equal order.
By Theorem \ref{arc-stabilizer},   $M_{\{u,v\}}$ is nonabelian. Thus  $(T_i)_{\{u,v\}}$ is nonabelian for all $i$. Then the lemma follows.
\qed

\vskip 5pt

For the case where $\Ga$ is a bipartite graph, we let $G^+$ be the subgroup of $G$ preserving the bipartition of $\Ga$.
Then $|G:G^+|=2$, and each bipartite half of $\Ga$ is a $G^+$-orbit on $V$.
\begin{lemma}\label{biquasi-case}
Assume  that the graph $\Ga=(V,E)$ is $(G,2)$-arc-transitive, and  $G$ has type {\rm PA} on $E$. Then either  $\Ga\cong \K_{d,d}$,
or one of the following holds:
\begin{enumerate}
\item[(1)] $G$ is quasiprimitive on $V$;
\item[(2)] $\Ga$ is bipartite, and $G^+$ is faithful and quasiprimitive on each bipartite half  of $\Ga$.
\end{enumerate}
\end{lemma}
\proof
Since $G$ is primitive on $E$, every minimal normal subgroup of $G$ is transitive on $E$, and so has at most two orbits on $V$.
If $\Ga$ is not bipartite  then quasiprimitive on $V$.

Now let  $\Ga$ be bipartite with bipartition say $V=V_1\cup V_2$. Note that $G_v\le G^+$ for each $v\in V$. Then
$G^+$ is locally-primitive on $\Ga$.
Suppose  that $\Ga\not\cong \K_{d,d}$. Then, by \cite{Prag-o'Nan-bi}, $G^+$ is faithful on both $V_1$ and $V_2$,
and either (2) of this lemma holds, or the unique minimal normal subgroup of $G$ is a direct product $M_1\times M_2$, where
 $M_1$ and $M_2$ are normal in $G^+$ and
 conjugate in $G$, and $M_i$ is intransitive on $V_i$ for $i=1,2$.
 For the latter case, if $M_1$ is intransitive on $V_2$ then
 $M_1$ is semiregular on $V$ by \cite[Lemma 5.1]{GLP}; if $M_1$ is transitive on $V_2$ then $M_2$  is semiregular on $V_2$.
 These two cases all contradict Lemma \ref{G-E-PA}. Thus $G^+$ is  quasiprimitive on both $V_1$ and $V_2$.
\qed

As permutation groups
on $V$ and on $E$, the types of $G$ (and $G^+$) have been determined in \cite{Giu-Li-edge-prim}. Then by Lemma \ref{biquasi-case} and combining with the reduction theorems for
2-arc-transitive graphs
given by Preager \cite{Prag-o'Nan, Prag-o'Nan-bi}, we get the following result.

\begin{lemma}\label{reduct-1}
Assume  that the graph $\Ga=(V,E)$ is $(G,2)$-arc-transitive.
Suppose that $\Ga\not\cong \K_{d,d}$. If $G$ is not almost simple, then $G$ has type {\rm PA} on $E$ and
either
\begin{enumerate}
 \item[(1)] $G$ is quasiprimitive and of type {\rm PA} on $V$; or
 \item[(2)] $\Ga$ is bipartite, $G^+$ is faithful and quasiprimitive on each bipartite half  of $\Ga$ with type {\rm PA}.
   \end{enumerate}
\end{lemma}

\vskip 10pt

Now we are ready to give a proof of Theorem~\ref{main-result}.
\begin{theorem}\label{main-1}
Let $\Ga=(V,E)$ be a connected $d$-regular graph for some $d\ge 3$, and let $G\le \Aut\Ga$.
 Assume  that $\Ga$ is both $G$-edge-primitive and $(G,2)$-arc-transitive. Then either $\Ga\cong \K_{d,d}$, or $G$ is almost simple.
\end{theorem}
\proof
Assume that  $\Ga\not\cong \K_{d,d}$, and let $\{u,v\}\in E$.
By the 2-arc-transitivity of $G$ on $\Ga$, we know that $G_v^{\Ga(v)}$ is a $2$-transitive permutation group of degree $d$.

Let $M=\soc(G)=T_1\times \cdots \times T_k$, where  $T_i$ are isomorphic nonabelian simple groups.
Then $M_v\unlhd G_v$, and $M_v\ne 1$ by Lemma \ref{nor-sub} or \ref{G-E-PA}.
Thus $M_v^{\Ga(v)}$ is a transitive normal subgroup of $G_v^{\Ga(v)}$.

\vskip 5pt

 Assume that $M_v^{\Ga(v)}$ is primitive on $\Ga(v)$. Noting that $G$ is transitive on $V$, we conclude that  $M_w^{\Ga(w)}$ is primitive for every $w\in V$.
Thus $\Ga$ is $M$-locally primitive. Then, by Lemma \ref{G-E-PA} and \cite[Lemma 5.1]{GLP}, we conclude that $k=1$, and so $G$ is almost simple.

\vskip 10pt

Next assume  that $M_v^{\Ga(v)}$ is imprimitive on $\Ga(v)$.

Note that every non-trivial normal subgroup of an almost simple $2$-transitive group  is primitive.
Then $G_v^{\Ga(v)}$ is an affine  $2$-transitive group, and by Lemma \ref{nor-sub-affine}, $M_v^{\Ga(v)}$ is a soluble Frobenius group and
$(M_v^{\Ga(v)})_u$ is cyclic. Set $(M_v^{\Ga(v)})_u\cong \ZZ_e$  and
$\soc(G_v^{\Ga(v)})\cong \ZZ_r^l$ for a prime $r$ and
integer $l\ge 1$ with $d=r^l$. Then $e$ is a divisor of $r^l-1$, and $e<r^l-1$.

Assume that $e=1$. Then  $M_v^{\Ga(v)}=\soc(G_v^{\Ga(v)})\cong \ZZ_r^l$, and so $M_v^{\Ga(v)}$ is regular on $\Ga(v)$. By \cite[Lemma 2.3]{567},  $M_v$ is faithful and hence regular on $\Ga(v)$, and thus   $M_{uv}=1$, which
contradicts Corollary \ref{M-v>d}. Thus $e\ne 1$.

Note that $e$ is a proper divisor of $d-1=r^l-1$. Neither $d$ nor $d-1$ is a prime, in particular, $l>1$ and $d=r^l\ge 9$. Thus $G_v^{\Ga(v)}$ has no   normal subgroup
isomorphic to a projective special linear group of dimension $\ge 2$.
By Theorem \ref{double-star},
$G_{uv}^{[1]}=1$, and so $M_{uv}^{[1]}=1$.

Let $x\in G_{\{u,v\}}\setminus G_{uv}$. Then
$(u,v)^x=(v,u)$, this implies that $M_v^{\Ga(v)}$ and $M_u^{\Ga(u)}$ are permutation isomorphic.
In particular, $(M_u^{\Ga(u)})_v\cong (M_v^{\Ga(v)})_u=\ZZ_e$.
Since $M_v^{[1]}\cap M_u^{[1]}=M_{uv}^{[1]}=1$, we know that
 $M_{uv}$ is isomorphic to  a subgroup of $(M_{uv}/M_u^{[1]})\times (M_{uv}/M_v^{[1]})$. Note that $M_{uv}/M_v^{[1]}\cong (M_v^{\Ga(v)})_u$ and
$M_{uv}/M_u^{[1]}\cong (M_u^{\Ga(u)})_v$. Then $M_{uv}$ is isomorphic to a subgroup of $\ZZ_e\times\ZZ_e$.
In particular, $M_{uv}$ is abelian. Then, by Theorem \ref{arc-stabilizer}, $M$ is transitive on the arcs of $\Ga$, and so $M_{\{u,v\}}=M_{uv}.2$.

If $e$ is a power of $2$ then, by  Theorem \ref{arc-stabilizer}, $M\cong \PSL(2,r^l)$, $\Ga\cong \K_{r^l+1}$; however, in this case, $M$ is locally primitive on $\Ga$, a contradiction. Thus $e$ has odd prime divisors. Let $s$ be an odd prime divisor of $e$, and $S$ be a Sylow $s$-subgroup of $M_{uv}$.  Then, noting that  $M_{\{u,v\}}=M_{uv}.2$, we know that
 $S$ is also a Sylow $s$-subgroup of $M$ by Theorem \ref{arc-stabilizer}. Thus $S=S_1\times \cdots\times S_k$, where $S_i$ is  a Sylow $s$-subgroup of $T_i$ for $1\le i\le k$. Since
$M_{uv}$ is isomorphic to a subgroup of $\ZZ_e\times\ZZ_e$, we know that
$M_{uv}$ has no subgroup isomorphic to $\ZZ_s^3$. It follows that $k\le 2$.

\vskip 5pt

Now we deduce a contradiction by supposing that $k=2$.

Let $k=2$.
Since $G\le (\Aut(T_1)\times \Aut(T_1)){:}2$, we have
\[G_{\{u,v\}}/M_{\{u,v\}}=G_{\{u,v\}}/(M\cap G_{\{u,v\}})\cong MG_{\{u,v\}}/M=G/M\le (\Out(T_1)\times \Out(T_1)){:}2.\]
It follows that $G_{\{u,v\}}/M_{\{u,v\}}$ is soluble, and so $G_{\{u,v\}}$ is soluble as $M_{\{u,v\}}$ is   soluble.
Thus $(G_v^{\Ga(v)})_u$ is soluble, and  $G_v^{\Ga(v)}=\soc(G_v^{\Ga(v)}){:}(G_v^{\Ga(v)})_u$ is also soluble.
Checking the soluble affine $2$-transitive groups, by Lemma \ref{nor-sub-affine},
 $(G_v^{\Ga(v)})_u\le \GammaL(1,r^l)$  or  $\ZZ_e\cong (M_v^{\Ga(v)})_u\le \Z((G_v^{\Ga(v)})_u)\cong \ZZ_2$.
 Note that $(M_v^{\Ga(v)})_u$  is a reducible subgroup of $(G_v^{\Ga(v)})_u$.
 Recalling that $e$ is not a power of $2$, the latter case does not occur.

Since $|M_{\{u,v\}}:M_{uv}|=2$, we have $M_{\{u,v\}}\not\le G_{uv}$,
and so $G_{uv}\ne M_{\{u,v\}}G_{uv}\le G_{\{u,v\}}$. Then $M_{\{u,v\}}G_{uv}= G_{\{u,v\}}$, and $G=MG_{\{u,v\}}=MG_{uv}$.
Recalling that $M=T_1\times T_2$, it follows that $G_{uv}$ acts transitively on $\{T_1,T_2\}$ by conjugation. Let $H$ be the kernel of this action. Then $|G_{uv}:H|=2$, and each $T_i$ is normalized by $H$. For $h\in H$,
\[((T_i)_v)^h=(T_i\cap G_v)^h=T_i^h\cap (G_v)^h=T_i\cap G_v=(T_i)_v, \, \,i=1,2.\]
This implies that $H$ normalizes each $(T_i)_v$. Then $(T_i)_v^{\Ga(v)}$ is normalized by $H^{\Ga(v)}$. Note that $(T_i)_v^{\Ga(v)}$ is a normal
subgroup of $M_v^{\Ga(v)}=\soc(G_v^{\Ga(v)}){:}(M_v^{\Ga(v)})_u$, and $e=|(M_v^{\Ga(v)})_u|$ is a proper divisor of $r^l-1$. Let $K_i$ be the Sylow $r$-subgroup of $(T_i)_v^{\Ga(v)}$. Then $K_i$ is normalized by $H^{\Ga(v)}$, of course, $K_i\le \soc(G_v^{\Ga(v)})$ and $K_1\cap K_2=1$.

Recalling that $|G_{uv}:H|=2$, we have $|(G_v^{\Ga(v)})_u:H^{\Ga(v)}|\le 2$.
Since $G_v^{\Ga(v)}$ is $2$-transitive, $|(G_v^{\Ga(v)})_u|$ is divisible by $r^l-1$, and so $|H^{\Ga(v)}|$ is divisible by ${r^l-1}\over 2$.
Note that ${{r^l-1}\over 2}>{r^l\over 2}-1\ge r^{l-1}-1$.
Then $|H^{\Ga(v)}|$ is not a divisor of $r^b-1$ for all  $ 1\le b<l$.
Then, by \cite[Lemma 5.1]{Hering}, $H^{\Ga(v)}$ is irreducible on $\soc(G_v^{\Ga(v)})$. It implies that $K_1=K_2=1$, and thus
$(T_i)_v^{\Ga(v)}\le (M_v^{\Ga(v)})_u$ for $i=1,2$. Let $u$ run over $\Ga(v)$.
It follows that $(T_i)_v^{\Ga(v)}=1$, and hence $(T_i)_v\le M_v^{[1]}$, $i=1,2$.
Since $M$ is transitive on $V$, by \cite[Lemma 2.3]{567}, we have $(T_1)_v=(T_2)_v=1$, which contradicts Lemma \ref{G-E-PA}. This completes the proof.
\qed

\vskip 5pt

As consequence of Theorems \ref{arc-stabilizer} and \ref{main-1},  an edge-primitive graph
of prime valency is $2$-arc-transitive, and then it has almost simple automorphism group if it is not a complete bipartite graph.
See also \cite{pval}.

\begin{corollary}\label{p-val}
Assume that $d$ is a prime and $\Ga\not\cong \K_{d,d}$. Then $G$ is almost simple, and  either $G=\PSL(2,d)$ with $d>11$ and $\Ga\cong \K_{d+1}$ or  $G$ is transitive on the $2$-arcs of $\Ga$.
\end{corollary}
\proof
 Note that $G$ is transitive on the arc set of $\Ga$. Let $\{u,v\}\in E$. By Theorem \ref{main-1}, it suffices to deal with the case where $G_v^{\Ga(v)}$ is not $2$-transitive.

Suppose that $G_v^{\Ga(v)}$ is not $2$-transitive. Then
$G_v^{\Ga(v)}\cong \ZZ_d{:}\ZZ_l$ with $l<d-1$ and $l$ a divisor of $d-1$. If $l=1$ then
$G_v\cong \ZZ_d$ by \cite[Lemma 2.3]{567}, and so $G_{uv}=1$, which contradicts Corollary \ref{M-v>d}. Then $l>1$, and so $d\ge 5$. By Theorem \ref{double-star},
$G_{uv}^{[1]}=1$. Then $G_{uv}$ is isomorphic to  a subgroup of $(G_{u}^{\Ga(u)})_v\times (G_{v}^{\Ga(v)})_u\cong \ZZ_l\times \ZZ_l$. Thus $G_{uv}$ is abelian. By Theorem \ref{arc-stabilizer},  $\Ga\cong \K_{d+1}$, $\soc(G)\cong \PSL(2,d)$, $\soc(G)_{v}\cong \ZZ_d{:}\ZZ_{d-1\over 2}$ and $\soc(G)_{\{u,v\}}\cong \D_{d-1}$. If $G\cong \PGL(2,d)$ then $G$ is transitive on the $2$-arcs of $\Ga$, which is not the case. Thus $G\cong \PSL(2,d)$, and so $d>11$ by the maximality of $G_{\{u,v\}}$.
\qed

\vskip 30pt

\section{Examples}
Let $\Ga=(V,E)$ be a connected $d$-regular graph, where $d\ge 3$. Let $v\in V$ and $G\le \Aut\Ga$.
Assume  that $\Ga$ is $(G,2)$-arc-transitive. Choose an  integer $s\ge 2$ such that   $\Ga$ is $(G,s)$-arc-transitive but not $(G,s+1)$-arc-transitive; in this case, we call $\Ga$ a $(G,s)$-transitive graph.
Then  $s\le 7$ by \cite{Weiss81b}.
If $G_v$ is faithful on $\Ga(v)$ then $s\le 3$
 by Theorem \ref{double-star}, and $s=3$ yields that $d=7$ and $G_v\cong \A_7$ or $\S_7$, see  \cite[Proposition 2.6]{Li-04}.
  This leads to the following interesting problem:
{\em Do there exist $3$-arc-transitive graphs with faithful stabilizers?}
We next answer this problem by giving several examples of edge-primitive graphs which are $3$-transitive and have faithful stabilizers.

\vskip 5pt

The first example is the Hoffman-Singleton graph, which has valency $7$, order $50$ and automorphism group $G=T.2$, where $T=\PSU(3,5)$. Let $X=T$ or $G$. For an edge ${\{u,v\}}$ of  this graph,  $X_v\cong \A_7$ or $\S_7$ and $X_{\{u,v\}}\cong \M_{10}$ or $\PGammaL(2,9)$, which are  maximal subgroups of  $X$.
Thus the  Hoffman-Singleton graph is both $X$-edge-primitive and $(X,2)$-arc-transitive.
To see the $3$-arc-transitivity, we fix an edge $\{u,v\}$ and consider the action of the arc-stabilizer $X_{uv}$ ($\cong \A_6$ or $\S_6$) on $\Ga(u)\cup \Ga(v)$. By the $2$-arc-transitivity of $X$, we have two faithful transitive  actions of $X_{uv}$ on $\Ga(u)$ and $\Ga(v)$, respectively. Let $v_1\in \Ga(v)\setminus\{u\}$ and $x\in X_{\{u,v\}}\setminus X_{uv}$. Then $u_1:=v_1^x\in \Ga(u)\setminus \{v\}$, and \[(X_{uv})_{u_1}=(X_{\{u,v\}})_{u_1}=(X_{\{u,v\}})_{v_1^x}
 =((X_{\{u,v\}})_{v_1})^x=((X_{uv})_{v_1})^x.\] By the choice of $x$, we know that $(X_{uv})_{v_1}$ and $((X_{uv})_{v_1})^x$ are not conjugate in $X_{uv}$, and so do for $(X_{uv})_{v_1}$ and $(X_{uv})_{u_1}$. This implies  that the actions  of $X_{uv}$ on $\Ga(u)$ and $\Ga(v)$ are not equivalent. Thus $(X_{uv})_{v_1}$ acts on $\Ga(u)\setminus \{v\}$ without fixed-points, this yields that $(X_{uv})_{v_1}$ is transitive on  $\Ga(u)\setminus \{v\}$. It follows that the Hoffman-Singleton graph is $(X,3)$-arc-transitive.

 \vskip 5pt

In general, combining with \cite[Proposition 2.6]{Li-04}, a similar argument as above yields the following result.

\begin{lemma}\label{3-arc-t-f-stab}
Let $\Ga=(V,E)$ be a connected $d$-regular graph for $d\ge 3$,   $\{u,v\}\in E$ and $G\le \Aut\Ga$.
If $\Ga$ is $(G,2)$-arc-transitive  and  $G_v$ is faithful on $\Ga(v)$, then $\Ga$ is $(G,3)$-arc-transitive if and only if $d=7$, $\soc(G_v)\cong \A_7$ and $G_{\{u,v\}}\not\cong \S_6$, i.e. $G_{\{u,v\}}\cong \PGL(2,9)$, $\M_{10}$ or $\Aut(\A_6)$.
\end{lemma}

We next give another example.

\begin{example}\label{ON}
{\rm By the information given in the Atlas \cite{Atlas} for the O'Nan simple group $\ON$,
there exactly two conjugacy classes $\mathcal{C}_1$ and $\mathcal{C}_2$ of (maximal) subgroups isomorphic to $\A_7$, which are merged into one class in $\ON.2$.
Further, there are $H\in \mathcal{C}_1$ and involutions $x_1,x_2\in \ON.2\setminus\ON$ such that $(H\cap H^{x_i}){:}\l x_i\r$ all are  maximal subgroups of $\ON.2$ with $(H\cap H^{x_1}){:}\l x_1\r\cong \PGL(2,9)$ and $(H\cap H^{x_2}){:}\l x_2\r\cong \PSL(2,7){:}2$. Define two bipartite graphs $\Ga_1=(V,E_1)$ and $\Ga_2=(V,E_2)$ with vertex set $V=\mathcal{C}_1\cup \mathcal{C}_2$ and edge sets
\begin{enumerate}
\item[] $E_1=\{\{H_1,H_2\}\mid H_1\in \mathcal{C}_1, H_2\in \mathcal{C}_2, H_1\cap H_2\cong \A_6\}$;
\item[] $E_2=\{\{H_1,H_2\}\mid H_1\in \mathcal{C}_1, H_2\in \mathcal{C}_2, H_1\cap H_2\cong \PSL(2,7)\}$.
\end{enumerate}
Then $\Ga_1$ and $\Ga_2$ are both $\ON.2$-edge-primitive and $(\ON.2,2)$-arc-transitive, which have valency $7$ and $15$ respectively.
By Lemma \ref{3-arc-t-f-stab}, only $\Ga_1$ is $(\ON.2,3)$-arc-transitive.}\qed
\end{example}

\begin{lemma}\label{auto-on}
Let $\Ga_1$ be as in Example {\rm \ref{ON}}. Then $\Aut\Ga_1=\ON.2$.
\end{lemma}
\proof
Let $G=\Aut\Ga_1$. Then $G\ge \ON.2$. By Theorem \ref{main-result}, $G$ is almost simple,  and
so $\ON\le \soc(G)$. Let $\{u,v\}$ be an edge of $\Ga_1$. Then $G_v^{\Ga(v)}\cong \A_7$ or $\S_7$, and $G_{uv}^{[1]}=1$ by Theorem \ref{double-star}. Thus, by the group extensions ($\circledast$) in Section \ref{pre}, we conclude that $|G_v|$ has no prime divisor other than $2$, $3$, $5$ and $7$. Since $\ON.2$ is transitive on the vertices of $\Ga_1$, we have $G=(\ON.2)G_v$. It follows that $|\ON|$ and $|\soc(G)|$ have the same prime divisors.
Using \cite[Corollary 5]{LPS-00}, we get $\soc(G)=\ON$, and so $G=\ON.2$.
\qed

\end{document}